\documentclass[12 pt]{amsart}
\usepackage{color}
\title{Heaviness in Toral Rotations}

\newtheorem{theorem}{Theorem}

\author{Yuval Peres}
\address{Yuval Peres\\ Microsoft Research\\ Microsoft Corporation\\ Redmond, WA}
\email{peres@microsoft.com}

\author{David Ralston}
\address{David Ralston\\ Department of Mathematics\\ The Ohio State University\\ 231 West 18th Ave.\\ Columbus, OH 43210}
\email{ralston@math.ohio-state.edu}


\date{\today}
\newcommand{\E}{\epsilon}
\newcommand{\A}{\alpha}
\newcommand{\T}{\mathbb{T}}
\newcommand{\R}{\mathbb{R}}
\newcommand{\Z}{\mathbb{Z}}
\newcommand{\N}{\mathbb{N}}
\newcommand{\Q}{\mathbb{Q}}
\newcommand{\lowerdim}{\underline{\dim}_M}
\newcommand{\upperdim}{\overline{\dim}_M}

\begin{document}

\begin{abstract}
We investigate the dimension of the set of points in the $d$-torus which have the property that their orbit under rotation by some $\A$ hits a fixed closed target $A$ more often than expected for all finite initial portions.  An upper bound for the lower Minkowski dimension of this \emph{strictly heavy set} $H(A,\A)$ is found in terms of the upper Minkowski dimension of $\partial A$, as well as $k$, the Diophantine approximability from below of the Lebesgue measure of $A$:
\[\lowerdim(H(A,\A)) \leq \upperdim{\partial A}+ \frac{d-\upperdim{\partial A}}{k}.\]  The proof extends to translations in compact abelian groups more generally than just the torus, most notably the $p$-adic integers.
\end{abstract}

\maketitle

\section{Introduction and Results}

The topic of \emph{heaviness} is broadly concerned with the maintenance of a fixed lower bound for a sequence of partial sums.  For example, consider $A$, a closed subset of the $d$-torus $\T^d = \R^d / \Z^d$, and let $\A \in \T^d$.  The transformation $T(x) = x+\A \mod{1}$ from the $d$-torus to itself preserves Lebesgue measure (denoted $\mu$), so we define the \emph{strictly heavy set} of $A$ with respect to rotation by $\A$:
\[H(A,\A) = \left\{ x \in \T^d : \sum_{i=0}^{n-1} \left( \chi_A(x + i \A \mod{1}) - \mu(A) \right) > 0, \hspace{.1 in} \forall n \in \N \right\}.\]
Heaviness in this context becomes a statement about the recurrence of the set $A$ to itself under the map $T$ in terms of Schnirelmann density.  In \cite{peres}\cite{ralston-houston}, it is shown that this set is nonempty if the inequality is not strict.  Clearly, then, if $\mu(A) \notin \Q$, then the set $H(A,\A)$ is not empty, as the change to a strict inequality is inconsequential in this case. 

We refer to \cite{falconer} for details regarding fractal dimension, and present only those definitions we will directly need.  For a set $S \subset \T^d$ define $\#(S)$ to be the cardinality of that set, and for $\E>0$, define  \[S_{\E} = \left\{x \in \T^d : \inf_{y \in S} \left|x-y\right| \leq \E \right\}.\]  Let $\mathcal{S}$ be the collection of all finite subsets of $S$ whose points are separated by at least $2\E$.  The \emph{packing number} of $S$ (given $\E$) is given by  \[N(\E) = \max\left\{ \#(S') : S' \in \mathcal{S}\right\}.\]  The \emph{upper/lower Minkowski dimension}  of $S$ are defined by 
\begin{align*}
\upperdim(S) &= \limsup_{\E \rightarrow 0} \left( \frac{\log N(\E) }{-\log \E}\right), \\
\lowerdim(S) &= \liminf_{\E \rightarrow 0} \left( \frac{\log N(\E)}{-\log \E}\right),
\end{align*}
and the \emph{upper/lower $s$-dimensional content} of $S$ are defined by
\begin{align*}
\overline{\mu_s}(S) &= \limsup_{ \E \rightarrow 0} \frac{\mu(S_{\E}) - \mu(S)}{\E^{d-s}},\\
\underline{\mu_s}(S)& = \liminf_{ \E \rightarrow 0} \frac{\mu(S_{\E})-\mu(S)}{\E^{d-s}}.
\end{align*}

If $\overline{\mu_s}(A)< \infty$, then it follows that for some fixed constant $c_1$, for all sufficiently small $\E$:
\begin{equation}
\label{dimension thing}
\mu \left( A_{\E} \right) - \mu(A) \leq c_1 \epsilon^{d-s}.
\end{equation}
Regardless of content, if $s=\upperdim(\partial A)$, then for any $\tau>0$, for all sufficiently small $\E$:
\begin{equation}
\label{dimension thing 2}
\mu \left( A_{\E} \right) - \mu(A) \leq \epsilon^{d-s-\tau}.
\end{equation}

A real number $\gamma$ is said to be \emph{approximable to order} $k$ if there is some constant $c_2$ and an infinite sequence of rational $p_i/q_i$ such that \[\left| \gamma - \frac{p_i}{q_i} \right| < \frac{c_2}{q_i^k}.\]  In this vein, we define $\gamma$ to be \emph{approximable from below to order} $k$ if there is an infinite sequence of rational $p_i/q_i$ such that \[0 \leq \gamma - \frac{p_i}{q_i} < \frac{c_2}{q_i^k}.\]  We note that every irrational number is approximable from below to order $2$ (an elementary result in the theory of continued fractions: see \cite[Ch. 2, \S 7]{khinchin}).

Our result is the following:
\begin{theorem}\label{maintheorem}
For any closed $A \subset \T^d$, where $\mu(A) \notin \Q$ is approximable from below to order $k$, for almost every $\A$,
\[\lowerdim(H(A,\A)) \leq \upperdim(\partial A) + \frac{d-\upperdim(\partial A)}{k}.\]  If $\mu(A) \in \Q$, then for almost every $\A$,
\[\lowerdim(H(A,\A)) \leq \upperdim{\partial A}.\]
\end{theorem}

We note that in the event that $\mu(A) = 1$, then as $A$ is closed, $A = \T^d$, $\partial{A}=\emptyset$, and $H(A,\A) = \emptyset$, satisfying the conclusion of the theorem.  If $\mu(A)=0$, then $H(A,\A)=A$, and the theorem is also trivially true.

\section{Proof of Theorem \ref{maintheorem}}

Our proof of Theorem \ref{maintheorem} will apply in greater generality than translations in the $d$-torus, so we fix $G$, an abelian group with metric $\left| \hspace{.1 in} \right|$, so that translations are isometries.  With respect to the topology induced by this metric, assume $G$ to be compact.  Let $\mu$ be the Haar measure on $G$, normalized so that $\mu(G)=1$.  Assume a certain regularity between the metric and the measure: there exists real $d$ and two finite, positive constants $c_3$ and $c_4$ so that for sufficiently small $\E$,
\begin{equation}
\label{regularity}
c_3 \E^d \leq \mu( \left\{x\right\}_{\E}) \leq c_4 \E^d.
\end{equation}
This condition on $\mu$ is called \emph{$d$-dimensional Ahlfors regularity} in \cite[p. 62]{david-semmes}.  Then the definitions of dimension and content transfer immediately, so we let $A \subset G$ be closed, with $\mu(A)$, $\mu(G \setminus A) \neq 0$ so that $\upperdim(\partial A) = \psi \geq d-1$. 

Define the variables $X_i: G \times G \rightarrow \R$ by \[X_i(x,g) = \chi_A(x+ig)-\mu(A),\] and define the set \[h_X(n,g)= \left\{x \in G : \sum_{i=0}^{j-1} X_i(x,g) > 0, \hspace{.1 in} j=1,2,\ldots,n\right\}.\]  Then we have \[H(A,g) = \bigcap_{n=1}^{\infty} h_X(n,g).\]

The following theorem may be found in \cite[p. 457]{loeve}:
\begin{theorem}
\label{loeve theorem}
{Suppose the random variables $V_i$, $i=0,1,2,\ldots$, are orthogonal and of mean zero.  Then: \[\mathbb{E}\left(\max_{k=1,2,\ldots,n} \left| \sum_{i=0}^{k-1} V_i \right| \right)^2 \leq \left(\frac{\log(4n)}{\log 2}\right)^2  \sum_{i=0}^{n-1} \mathbb{E}\left| V_i\right|^2.\]}
\end{theorem}

Assume that $\mu(A) \notin \Q$ is approximable from below to order $k$, where $\left\{p_i/q_i\right\}_{i=1}^{\infty}$ is an increasing sequence such that \[0 \leq \mu(A) - \frac{p_i}{q_i} < \frac{c_2}{q_i^k}.\]  The case $\mu(A) \in \Q$ is simpler, and will be addressed later.  With this sequence fixed, set \[n_i = q_i^{\frac{2k}{d-\psi}}, \hspace{.1 in} \E_i = n_i^{-\frac{1}{2}}.\]
Define the random variables \[Y_j^{(i)}(x,g) = \chi_{A_{\E_i}}(x+jg) - \frac{p_i}{q_i},\] and the sets \[h_Y(n,g) = \left\{x \in G : \sum_{i=0}^{j-1} Y_i(x,g)>0, \hspace{.1 in} j=1,2,\ldots,n\right\}.\]

Note that different superscripts $i$ substantially change the variables $Y$ by changing both $\epsilon_i$ and $p_i/q_i$.  Also note that as $p_i/q_i< \mu(A)$, if $x \in h_X(n_i,g)$, then $\{x\}_{\E_i} \subset h_Y(n_i,g)$.  Recall that for any $i$ we have $H(A,g) \subset h_X(n_i,g)$.  Therefore, if $N(\E_i)$ is the packing number of $H(A,g)$,
\begin{equation}
\label{lower bound}
N(\E_i) c_3 \E_i^d \leq \mu\left( h_Y(n_i,g) \right).
\end{equation}

Assume that the upper $\psi$-dimensional content of $A$ is finite, and define \[Z_j^{(i)}(x,g) = \chi_{A_{\E_i}}(x+jg) - \mu \left( A_{\E_i}\right).\] The variables $Z_j$ are orthogonal over their domain $G \times G$, so by applying a H\"{o}lder inequality to Theorem \ref{loeve theorem}, we obtain:

\[\mathbb{E} \left( \max_{k=1,2,\ldots,n_i} \left| \sum_{j=0}^{k-1} Z_j^{(i)}(x,g) \right| \right) \leq c_4 \log(n_i) \sqrt{n_i},\] where $c_5$ is a constant which is independent of $i$ (the variance of our variables $Z_j^{(i)}$ may always be bounded by $1/4$, regardless of $A$ and $i$).  Note that

\begin{align*}
Y_j^{(i)}(x,g) &= \chi_{A_{\E_i}}(x+jg) - \frac{p_i}{q_i},\\
&=  Z_j^{(i)}(x,g) + \left( \mu \left( A_{\E_i} \right) - \mu(A)\right) + \left( \mu(A) - \frac{p_i}{q_i}\right).
\end{align*}
So by integrating and applying Theorem \ref{loeve theorem}, Eq. \eqref{dimension thing}, and our approximability condition:
\begin{align*}
\mathbb{E}\left(\max_{k=1,2,\ldots n_i} \left| \sum_{j=0}^{k-1} Y_j^{(i)}(x,g) \right| \right)&< c_5 \log( n_i) n_i^{\frac{1}{2}} + c_3 n_i \E_i^{d-\psi} + c_2 \frac{n_i}{q_i^k},\\
&< c_6 n_i^{1-\frac{d-\psi}{2}}.
\end{align*}
for some constant $c_6$, where the last inequality follows from the relationship between $q_i$, $n_i$, and $\E_i$, and yields
\[\mathbb{E}\left(\max_{k=1,2,\ldots,n_i}\sum_{j=0}^{k-1}Y_j^{(i)}(x,g) - \min_{k=1,2,\ldots,n_i}\sum_{j=0}^{k-1}Y_j^{(i)}(x,g)\right) < 2c_6 n_i^{1-\frac{d-\psi}{2}}.\]

As the values \[S_k(x,g)=\sum_{j=0}^{k-1}Y_j^{(i)}(x,g)\] are all expressible as rational numbers with denominator $q_i$, the distinct values taken by $S_k(x,g)$ are all separated by at least $1/q_i$, so that:
\begin{align*}
\mathbb{E}\left(\# \left\{ S_j(x,g) \right\}_{j=1}^{n_i} \right)&\leq 2 c_6 q_i  n_i ^{1-\frac{d-\psi}{2}} \\
&= 2c_6 n_i^{1-\frac{(d-\psi)(k-1)}{2k}}.
\end{align*}

Consider now that for any $x \in G$, if for some $m < n_i$ we have $x+mg \in h_Y(n_i,g)$, then then value $S_m(x,g)$ is not repeated by any other $S_j(x,g)$ where $m < j \leq n_i$.  Define \[J_i(x,g) = \# \left\{j : x+jg \in h_Y(n_i,g), \hspace{.1 in} j=0,1,\ldots,n_i-1\right\},\] so that
\[J_i(x,g) \leq \# \left\{S_j(x,g) \right\}_{j=1}^{n_i}.\]  It is apparent, however, that \[\int_{G} J_i(x,g) d\mu(x) = n_i \mu \left( h_Y(n_i,g)\right).\]  By integrating against $g$, we altogether have that
\[\int_{G} \mu \left( h_Y(n_i,g)\right)d\mu(g) \leq 2 c_6 n_i ^{-\frac{(d-\psi)(k-1)}{2k}},\]
from which it follows via Fatou's Lemma that
\[
\int_{G}\left( \liminf_{i \rightarrow \infty} \mu\left(h_Y (A,n_i)\right) n_i^{\frac{(d-\psi)(k-1)}{2k}}\right)d\mu(g) < \infty
.\]
Therefore, for almost every $g \in G$, there is some constant $c_7$ depending only on $g$ so that for all sufficiently large $i$:
\begin{equation}\label{upper bound}
\mu\left(h_Y(n_i,g)\right) \leq c_7 n_i ^{-\frac{(d-\psi)(k-1)}{2k}}
\end{equation}

By combining Equations \eqref{lower bound} and \eqref{upper bound} we may finally estimate $\lowerdim(H(A,g))$ for this full-measure set of $g$:
\begin{align*}
\lowerdim(H(A,g)) &= \liminf_{\E \rightarrow 0} \frac{\log\left( N(\E) \right)}{- \log (\E)}\\
& \leq \liminf_{i \rightarrow \infty} \frac{\log \left( N(\E_i)\right)}{-\log (\E_i)}\\
& \leq \liminf_{i \rightarrow \infty} \frac{\log \left( c_8 n_i^{\frac{d}{2} - \frac{(d-\psi)(k-1)}{2k}}\right)}{-\log \left(n_i^{-\frac{1}{2}} \right)}\\
&= \psi + \frac{d-\psi}{k},
\end{align*} where $c_8$ is again some constant independent of $i$.

We remark that the $\left(\psi + \frac{d-\psi}{k}\right)$-dimensional content of $H(A,g)$ is finite (under the assumption that the $\upperdim(\partial A)$-dimensional content of $A$ was finite): note that $H(A,g)_{\E_i} \subset h_Y(n_i,g)$, while $H(A,g)$ is a null set for those $g$ which define $\mu$-ergodic translations (an immediate consequence of the Birkhoff Ergodic theorem).  By applying Equation \eqref{upper bound}, we see that the content is finite (possibly zero).

To complete the proof for the case $\mu(A) \notin \Q$, we need only address the case where the $\psi$-dimensional content of $A$ is infinite.  We cannot apply Equation \eqref{dimension thing} in this case, but we may carry out the same proof using Equation \eqref{dimension thing 2} in its place to show that \[\lowerdim(H(A,g)) \leq \psi + \frac{d-\psi}{k} + \tau\] for arbitrary $\tau>0$, for the same result with no information about dimensional content.

If $\mu(A) = p/q$, we do not need to involve the approximating term $\mu(A) - p_i/q_i$ in our sums; setting $n_i=i$, $\E_i=n_i^{-1/2}$, the values taken by the sums \[\sum_{j=0}^{k-1} \left( \chi_{A_{\E}} (x+jg) - \mu(A) \right)\] are already discrete.  So, all of the previous arguments are simplified, and the simplified expression \[\lowerdim{H(A,g)} \leq \upperdim{(\partial A)}\] is achieved for almost every $g$.  The statements regarding the content of $H(A,g)$ also follows exactly as in the case of irrational measure.

\section{Concluding Remarks}

In the simplest case $G=\R/\Z=[0,1)$ under rotation by irrational $\A$, with $A=[0,1/2]$, explicit techniques for finding the unique strictly heavy point and the structure of the set of non-strictly heavy points are developed in \cite{ralston-rotations} and applied to recurrence properties of cylinder flows in \cite{chaika-ralston}.  It follows from Theorem \ref{maintheorem} that the set of strictly heavy points is of finite zero-dimensional capacity, i.e. finite; that the set is actually a singleton depends on specific properties of rotations on the circle.

In $\T^d$, Equation \eqref{regularity} may be simplified to equality on both sides for a single constant.  However, the flexibility afforded by using different constants for the scaling of $\E$-balls allows us to apply our results to more general compact Abelian groups, notably the $p$-adic integers.  If information about the dimensional content of the heavy set is not required, however, we may relax the Ahlfors regularity of $\mu$, instead assuming that there are two functions $f$, $g$ so that \[f(\E) \E^d \leq \mu( \{x\}_{\E}) \leq g(\E)\E^d,\] and \[\lim_{\E \rightarrow 0} \frac{\log (f(\E))}{\log(\E)} = \lim_{\E \rightarrow 0} \frac{\log(g(\E))}{\log(\E)} = 0\] and still arrive at the same upper bound on the dimension of the heavy set.  

We have also assumed that translations are isometries.  This assumption may be relaxed to simply assuming that all translations are continuous with respect to the topology induced by some given metric.  In this case, compactness of $G$ gives that the family of translations is an equicontinuous family, so that there is some constant $c_9$ for which we have for all $x,g \in G$:
\[\left\{ x+g \right\}_{c_9^{-1} \E} \subset \left\{y+g: y \in \left\{x\right\}_{\E}\right\} \subset \left\{x+g\right\}_{c_9 \E}.\]

In this case, we need to adjust our definition of the variables $Y_i$ and $Z_i$ to be based off not the $\E$-neighborhood of $A$, but the $c_9 \E$-neighborhood of $A$.  Similarly to the statement about scaling measure of $\E$-balls, if we are not interested in statements about dimensional content, we need only know that the distortion of translates of $\E$-balls is negligible compared to $\log(\E)$.  The case of $G=\T^d$ is the most natural to consider initially, but each of these relaxations allows Theorem \ref{maintheorem} to apply in slightly more general contexts beyond the $d$-torus ($p$-adics, solenoids, etc.).

\section*{Acknowledgements}

The authors would like to thank Jon Chaika and Tim Austin for helpful comments regarding the clarity of exposition and for providing valuable references.  The second author wishes to thank Manfred Einsiedler and Michael Geline for helpful conversations.

\end{document}